\documentclass[12pt,leqno]{amsart}
\usepackage{amsmath,amsthm,amsfonts,amsbsy,multicol,fancybox}
\usepackage[utf8]{inputenc}
\usepackage[T1]{fontenc}

\usepackage{caption}\usepackage{subcaption}
\usepackage{graphicx}\usepackage{amssymb,latexsym}
\usepackage{epstopdf}
\usepackage{color}
\DeclareGraphicsExtensions{.pdf, .jpg, .pnp, .bmp, .fig}
\usepackage{tikz}
\usepackage[all]{xy}
\input xy
\xyoption{all}

\usepackage[hidelinks]{hyperref}

\usepackage{thmtools}

\newcounter{argument}
\newenvironment{argument}[1][\medskip]{%
\refstepcounter{argument}
\par\medskip
\noindent\phantomsection
\textbf{#1~\thesection.\arabic{argument}\,\,}\rmfamily\em}{\hspace{\fill}$\Box$\par\smallskip\noindent}

\newcommand{\bass}{\begin{argument}[Assumption]}\newcommand{\eass}{\end{argument}}
\newcommand{\bth}{\begin{argument}[Theorem]} \newcommand{\ethe}{\end{argument}}
\newcommand{\bre}{\begin{argument}[Remark]}      \newcommand{\ere}{\end{argument}}
\newcommand{\ble}{\begin{argument}[Lemma]}       \newcommand{\ele}{\end{argument}}
\newcommand{\bde}{\begin{argument}[Definition]}   \newcommand{\ede}{\end{argument}}
\newcommand{\bco}{\begin{argument}[Corollary]}     \newcommand{\eco}{\end{argument}}
\newcommand{\bpr}{\begin{argument}[Proposition]}  \newcommand{\epr}{\end{argument}}
\newcommand{\bexam}{\begin{argument}[Example]}\newcommand{\eexam}{\end{argument}}

\newcommand{\bpf}{\begin{proof}}\newcommand{\epf}{\end{proof}}
\newcommand{\barr}{\begin{array}}\newcommand{\earr}{\end{array}}
\newcommand{\beao}{\begin{eqnarray*}}\newcommand{\eeao}{\end{eqnarray*}\noindent}
\newcommand{\beam}{\begin{eqnarray}}\newcommand{\eeam}{\end{eqnarray}\noindent}
\newcommand{\beqq}{\begin{equation}}\newcommand{\eeqq}{\end{equation}\noindent}

\newcommand{\wt}{\widetilde}

\newcommand{\Rto}{R\to\infty}

\newcommand{\D}{\Delta}
  \newcommand{\ep}{\epsilon}

\newcommand{\lam}{\lambda}

\newcommand{\bfE}{{\mathbb E}}\newcommand{\bbE}{{\mathcal E}} 
\newcommand{\bbf}{{\mathcal F}}

\newcommand{\bbi}{{\mathbb I}}

\newcommand{\bbl}{{\mathcal L}}

 \newcommand{\bbN}{{\mathbb N}}

\newcommand{\bfP}{{\mathbb P}}
                              
 \newcommand{\bbR}{{\mathbb R}}

\graphicspath{{./../Images/}}

\begin{document}
\title[Numerical approximation for non-colliding particle systems]{Numerical approximation for non-colliding particle systems}

\author[I. S. Stamatiou]{I. S. Stamatiou}
\email{joniou@gmail.com, ioannis.stamatiou@ouc.ac.cy}

\begin{abstract}
We apply the semi-discrete method, c.f. \emph{N. Halidias and I.S. Stamatiou (2016), On the numerical solution of some non-linear stochastic differential equations using the semi-discrete method, Computational Methods in Applied Mathematics, 16(1)}, to a class of non-colliding particle systems. The proposed numerical scheme preserves the non-colliding property and strongly converges to the exact solution.
\end{abstract}

\date\today

\keywords{Explicit Numerical Scheme; Semi-Discrete Method; non-linear SDEs; Stochastic Differential Equations; Boundary Preserving Numerical Algorithm; Dyson Brownian motion; non-colliding particle system
 \newline{\bf AMS subject classification 2010:}  41A25; 60H10, 60H35, 60J60, 65C20, 65C30, 65J15, 65L20.}
\maketitle

\section{Introduction}\label{NCS:sec:intro}
\setcounter{equation}{0}

We are interested in the following system of stochastic differential equations (SDEs),
\beqq  \label{NCS-eq:systemSDEs}
X_t^{(i)} = X_0^{(i)} + \int_0^t  \left(\sum_{i\neq j}\frac{\gamma_{i,j}}{X_s^{(i)}-X_s^{(j)}} + b^{i}(X_s^{(i)})\right)ds + \sum_{j=1}^d\int_0^t \sigma_{i,j} dW_s^{(j)}, \quad i=1,\ldots,d,
\eeqq
where $X_0=(X_0^{(1)},\ldots,X_0^{(d)})^T\in \D_d=\{{\bf{x}}=(x^{(1)},x^{(2)},\ldots,x^{(d)})^T\in\bbR^d: x^{(1)}<x^{(2)}<\ldots<x^{(d)}\}$ almost surely (a.s.) and
$\{W_{t}\}_{t\geq0}$ is a $d$-dimensional Wiener process adapted to the filtration  $\{\bbf_t\}_{t\geq0}.$ The constants $\gamma_{i,j}$ satisfy $\gamma_{i,j}=\gamma_{j,i}\geq0$ with $\gamma_{i,i+1}>0$ for $i=1,\ldots,d-1$ and the functions $b^{i}(\cdot)$ are globally Lipschitz continuous or non-increasing with the property  $b^{i}(z)\leq b^{i+1}(z)$ for all $z\in\bbR$ and $\sigma_{i,j}$ are finite constants such that $\sigma^2:=\sup_{1\leq i\leq d}\sum_{k=1}^d\sigma_{i,k}^2\leq 2\gamma_{i,j}$. Under the above assumptions the system (\ref{NCS-eq:systemSDEs}) has a unique strong solution in $\D_d,$ c.f. \cite{graczyk_malecki:2014}, \cite{ngo_taguchi:2018}. We want to reproduce the non-colliding property of (\ref{NCS-eq:systemSDEs}). We use a fixed-time step explicit numerical method $(Y_{n}^{i+1,i})_{n\in\bbN}$, namely the semi-discrete method, to approximate the difference $X_{n+1}^{i+1,i}:=X_{t_{n+1}}^{(i+1)}-X_{t_{n+1}}^{(i)},$ which reads
\beqq  \label{NCS-eq:SDmethod}
Y_{n+1}^{i+1,i} =\sqrt{(\wt{Y}_{n+1}^{i+1,i})^2 + 4\gamma_{i+1,i}\D}, \quad n\in\bbN,
\eeqq
where 
$$
\wt{Y}_{n+1}^{i+1,i} =  e^{\alpha^i(\wt{Y}_n)\D} \left(\wt{Y}_n^{i+1,i} -\frac{\beta^i(\wt{Y}_n)}{\alpha^i(\wt{Y}_n)}(1-e^{-\alpha^i(\wt{Y}_n)\D})
+ \sqrt{\sum_{k=1}^d(\sigma_{i+1,k}-\sigma_{i,k})^2}e^{-\alpha^i(\wt{Y}_n)\D}\D B^{i+1}_n \right)
$$
and 
$$
\alpha^i(Y_n):=\frac{b^{i+1}(Y_n^{(i+1)})-b^{i}(Y_n^{(i)})}{Y_n^{(i+1)}-Y_n^{(i)}}- \sum_{k\neq i,i+1}\frac{\gamma_{i+1,i}}{Y_n^{i+1,k}Y_n^{k,i}}, \quad  \beta^i(Y_n):=\sum_{k\neq i,i+1}\frac{\gamma_{i,k}-\gamma_{i+1,k}}{Y_n^{i+1,k}},$$
with $Y_0=X_0; \D=t_{n+1}-t_{n}$ is the time step-size and $\D B^{i+1}:= B^{i+1}_{t_{n+1}}- B^{i+1}_{t_{n}}$ are the increments of a Wiener process. For the derivation of (\ref{NCS-eq:SDmethod}) see Section \ref{NCS:sec:proofs}.

Our main goal is to provide an explicit numerical method for the approximation of the solution of system (\ref{NCS-eq:systemSDEs}), which preserves the non-colliding property. The class of systems of SDEs of type (\ref{NCS-eq:systemSDEs}) has been studied thoroughly, see among others in \cite{Cepa1997} \cite{ramanan2018},  due to its practical importance: in the context of particle systems the term $\sum_{i\neq j}\frac{\gamma_{i,j}}{X^{(i)}-X^{(j)}}$ describes the repulsive force with which each of the $j$-th particle (located at $X^{(j)}$) acts on the $i$-th particle (located at $X^{(i)}$);
still few results on their numerical approximations are available: the explicit tamed Euler scheme  proposed in \cite{Li_Menon:2013} for the case of the Dyson Brownian motion, c.f. \cite{dyson:1962} ($\gamma_{i,j}=\gamma, \sigma_{i,j}=\delta_{i,j}\sigma_i $ where $\delta_{i,j}$ is the Dirac delta function) fails to preserve the non-colliding property; the only method we know of which preserves the non-colliding property is a semi-implicit EM method recently proposed by \cite{ngo_taguchi:2018}. 
 
The non-colliding property of (\ref{NCS-eq:systemSDEs}) is rephrased in the positivity property of the process $(X_{t}^{i+1,i});$ our method preserves this property by construction.

Moreover, we study the strong rate of convergence in $\bbl^p$-norm of the proposed method.

The case of numerical approximations of scalar SDEs with boundaries (like CIR process, mean-reverting CEV process, Wright-Fisher model) has been studied by many researchers, c.f. \cite{halidias_stamatiou:2015}, \cite{stamatiou:2018}, \cite{neuenkirch_szpruch:2014}, \cite{Dereich_et_al:2011}, \cite{alfonsi:2013}, \cite{moro_schurz:2007}. 

The proposed fixed-step method is explicit, strongly convergent, non-explosive and  preserves the non-colliding property. The semi-discrete method was originally proposed in \cite{halidias:2012} and further investigated  in  \cite{halidias_stamatiou:2016}, \cite{halidias:2014} for one-dimensional SDEs and  \cite{halidias:2015}, \cite{halidias:2015d} in the multivariate case. The basic ingredient of the semi-discrete method is the following: we freeze on each subinterval appropriate parts of the drift and diffusion coefficients of the solution at the beginning of the subinterval in order to obtain explicitly solved SDEs. Apparently the way of freezing (discretization) is not unique. Here, we freeze the nonlinear parts obtaining a linear SDE with explicit solution.
 
The outline of the article is the following. In Section \ref{NCS:sec:main} we present our main results, that is Theorem \ref{NCS-theorem:convergence} and Corollary \ref{NCS-corollary:convergence}, the proof of which are deferred to Section \ref{NCS:sec:proofs}. 

\section{Main results}\label{NCS:sec:main}

Let us restate the assumptions of the model (\ref{NCS-eq:systemSDEs}).
\bass\label{NSF:assA} 
\begin{itemize}
\item [i] $X_0\in \D_d$ a.s.
\item [ii] $\gamma_{i,j}=\gamma_{j,i}\geq0$ with $\gamma_{i,i+1}>0$ for $i=1,\ldots,d-1.$
\item [iii] $b^{i}(\cdot)$ are globally Lipschitz continuous with the property  $b^{i}(z)\leq b^{i+1}(z)$ for all $z\in\bbR.$
\item [iv] $\sigma_{i,j}$ are finite constants such that $\sigma^2:=\sup_{1\leq i\leq d}\sum_{k=1}^d\sigma_{i,k}^2\leq 2\gamma_{i,j}.$
\end{itemize}
\eass

Consider the process $X_t^{i+1,i}:=X_{t}^{(i+1)}-X_{t}^{(i)}, i=0,\ldots,d-1,$ where $X_{t}^{(0)}\equiv0,$ which satisfies the following SDE
\beam 
\nonumber
X_t^{i+1,i} &=& X_0^{i+1,i} + \int_0^t  \left(\frac{2\gamma_{i+1,i}}{X_s^{i+1,i}} - 
\sum_{k\neq i,i+1}\frac{\gamma_{i,k}X_s^{i+1,k}-\gamma_{i+1,k}X_s^{i,k}}{X_s^{i+1,k}X_s^{i,k}} + b^{i+1}(X_s^{(i+1)})-b^{i}(X_s^{(i)})\right)ds\\
\label{NCS-eq:system_diff_SDEs}
&&+ \sum_{j=1}^d\int_0^t (\sigma_{i+1,j}-\sigma_{i,j}) dW_s^{(j)}.
\eeam

The non-colliding property implies that $(X_t^{i+1,i})$ is positive.  In order to find the solution process, we use a splitting technique to get first, see also \cite{moro_schurz:2007},
\beam 
\nonumber
\wt{X}_t^{i+1,i} &=& \wt{X}_0^{i+1,i} + \int_0^t  \left( - \sum_{k\neq i,i+1}\frac{\gamma_{i,k}\wt{X}_s^{i+1,k}-\gamma_{i+1,k}\wt{X}_s^{i,k}}{\wt{X}_s^{i+1,k}\wt{X}_s^{i,k}} + b^{i+1}(\wt{X}_s^{(i+1)})-b^{i}(\wt{X}_s^{(i)})\right)ds\\
\label{NCS-eq:system_diff_split1}
&&+ \sum_{j=1}^d\int_0^t (\sigma_{i+1,j}-\sigma_{i,j}) dW_s^{(j)},
\eeam
and afterwards 
\beqq\label{NCS-eq:system_diff_split2}
X_t^{i+1,i}=\wt{X}_t^{i+1,i} + \int_0^t\frac{2\gamma_{i+1,i}}{X_s^{i+1,i}}ds = \sqrt{(\wt{X}_t^{i+1,i})^2 + 4\gamma_{i+1,i}t}.
\eeqq

Therefore instead of approximating directly the solution of (\ref{NCS-eq:system_diff_SDEs}) we work with SDE (\ref{NCS-eq:system_diff_split1}). Let $\wt{Y}_{t}^{i+1,i}$ be a fixed-step numerical approximation of $\wt{X}_{t}^{i+1,i}$ with time step size $\D$  and $\wt{Y}_{0}=\wt{X}_{0}$ such that  $\bfE\sup_{0\leq t\leq T} |\wt{Y}_{t}^{i+1,i} -\wt{X}_{t}^{i+1,i}|^2\leq C\D^{2p}$ where $p$ is the order of strong convergence. Then 
\beam\nonumber
|Y_{t}^{i+1,i} -X_{t}^{i+1,i}|^2&=&\left| \sqrt{(\wt{Y}_t^{i+1,i})^2 + 4\gamma_{i+1,i}t} - \sqrt{(\wt{X}_t^{i+1,i})^2 + 4\gamma_{i+1,i}t}\right|^2 \\
\nonumber &=&\left|(\wt{Y}_t^{i+1,i})^2 - (\wt{X}_t^{i+1,i})^2\right|\\
 \label{NCS-eq:SDsplit_part} &\leq&(|\wt{Y}_t^{i+1,i}|+ |\wt{X}_t^{i+1,i}|)|\wt{Y}_{t}^{i+1,i} -\wt{X}_{t}^{i+1,i}|, 
\eeam
which implies  $\bfE\sup_{0\leq t\leq T} |Y_{t}^{i+1,i} -X_{t}^{i+1,i}|^2\leq C\D^{p},$ provided that the first moments of $(\wt{Y}_{t}^{i+1,i})$ and $(\wt{X}_{t}^{i+1,i})$ are bounded.  Now, we discuss about the approximation $(\wt{Y}_{t}^{i+1,i})$ of $(\wt{X}_{t}^{i+1,i}).$ We propose an application of the semi-discrete method in the following way.

Let the equidistant partition $0=t_0<t_1<\ldots<t_N=T$ with step size $\D=T/N$ and consider the following process
\beam\nonumber
\wt{Y}_{t}^{i+1,i} &=& \wt{Y}_{0}^{i+1,i} +  \int_0^t\left(\frac{b^{i+1}(\wt{Y}_{\hat{s}}^{(i+1)})-b^{i}(\wt{Y}_{\hat{s}}^{(i)})}{\wt{Y}_{\hat{s}}^{i+1,i}} -  \sum_{k\neq i,i+1}\frac{\gamma_{i,k}}{\wt{Y}_{\hat{s}}^{i+1,k}\wt{Y}_{\hat{s}}^{i,k}}\right)\wt{Y}_{s}^{i+1,i}ds\\
\nonumber&&+ \int_0^t  \left( - \sum_{k\neq i,i+1}\frac{\gamma_{i,k} - \gamma_{i+1,k}}{\wt{Y}_{\hat{s}}^{i+1,k}}\right)ds + \sum_{j=1}^d\int_0^t (\sigma_{i+1,j}-\sigma_{i,j}) dW_s^{(j)}\\
 \label{NCS-eq:SDprocess}&=& \wt{Y}_{0}^{i+1,i} +  \int_0^t \left(\alpha^i(\wt{Y}_{\hat{s}})\wt{Y}_{s}^{i+1,i} - \beta^i(\wt{Y}_{\hat{s}}) \right)ds  + c^i \int_0^t  dB_s^{(i+1)},
\eeam
where $B_t^{(i+1)}:=\sum_{j=1}^d\int_0^t \frac{(\sigma_{i+1,j}-\sigma_{i,j})}{\sqrt{\sum_{k=1}^d(\sigma_{i+1,k}-\sigma_{i,k})^2}}dW_s^{(j)}$ is a new Wiener process,
\beam\label{NCS-eq:SDprocess_a}
\alpha^i (Y) &=&  \frac{b^{i+1}(Y^{(i+1)})-b^{i}(Y^{(i)})}{Y^{i+1,i}} -  \sum_{k\neq i,i+1}\frac{\gamma_{i,k}}{Y^{i+1,k}Y^{i,k}},\\
\label{NCS-eq:SDprocess_b}\beta^i (Y) &=& \sum_{k\neq i,i+1}\frac{\gamma_{i,k} - \gamma_{i+1,k}}{Y^{i+1,k}},\\ 
\label{NCS-eq:SDprocess_c} c^i  &=& \sqrt{\sum_{k=1}^d(\sigma_{i+1,k}-\sigma_{i,k})^2}
\eeam
and $\wt{Y}_{0}^{i+1,i}=\wt{X}_{0}^{i+1,i}$ a.s. Here $\hat{s}=t_n$ when  $s\in[t_n,t_{n+1})$ and denotes the freezing times; in particular we keep the diffusion part of (\ref{NCS-eq:system_diff_split1}) the same and freeze the drift in such a way that the produced SDE (\ref{NCS-eq:SDprocess}) is linear in the narrow sense with additive noise and unique strong solution given by (c.f. \cite[Sec 4.4.]{kloeden_platen:1995})
\beqq\label{NCS-eq:SDprocess_sol}
\wt{Y}_{t}^{i+1,i} = e^{\alpha^i(\wt{Y}_{\hat{s}}) t }\left(  \wt{Y}_{0}^{i+1,i} - \frac{\beta^i(\wt{Y}_{\hat{s}})}{\alpha^i(\wt{Y}_{\hat{s}})}(1-e^{-\alpha^i(\wt{Y}_{\hat{s}})t}) + c^i\int_0^t e^{-\alpha^i(\wt{Y}_{\hat{s}})s}dB_s^{(i+1)}\right). 
\eeqq
Note that the drift of (\ref{NCS-eq:SDprocess}) equals the drift of (\ref{NCS-eq:system_diff_split1}) for $\hat{s}=s.$
Our main result is the following.
\bth[Strong convergence of $\wt{Y}_t^{i+1,i}$ to $\wt{X}_t^{i+1,i}$]\label{NCS-theorem:convergence}
Let Assumption \ref{NSF:assA} hold, $(c^i)^2\geq(d-1)\sup_{i\neq j}\gamma_{i,j}$ and $\bfE\sup_{i\leq d-1}(|\wt{X}_{0}^{i+1,i}|^k\vee|\wt{X}_{0}^{i+1,i}|^{-k}\vee|\wt{Y}_{0}^{i+1,i}|^k\vee|\wt{Y}_{0}^{i+1,i}|^{-k})<A,$ for some $k\geq6$ and let $\bfE e^{C(\wt{X}_0^{i+1,i})^{-8}}<A_{X_0}$ for any $C>0$ where $A_{X_0}$ is a finite  constant. The semi-discrete scheme (\ref{NCS-eq:SDprocess}) converges strongly in the mean-square sense to the true solution of (\ref{NCS-eq:system_diff_split1}) with order of convergence $1/2$, that is 
$$
\bfE\sup_{0\leq t\leq T} |\wt{Y}_{t}^{i+1,i} -\wt{X}_{t}^{i+1,i}|^2\leq C\D,
$$ 
where $C$ is a constant independent of $\D.$
\ethe

The semi-discrete numerical scheme for the approximation of (\ref{NCS-eq:systemSDEs}) is the $d$-dimensional vector $Y_{t}$  where $Y_t^{(i)}=\sum_{j=0}^i Y_{t}^{j+1,j}$ with 
$Y_{t}^{1,0}=Y_{t}^{(1)},$
\beqq  \label{NCS-eq:SDmethod_t}
Y_{t}^{i+1,i} =\sqrt{(\wt{Y}_{t}^{i+1,i})^2 + 4\gamma_{i+1,i}t}, \quad n\in\bbN,
\eeqq
where $\wt{Y}_{t}^{i+1,i}$ is given in (\ref{NCS-eq:SDprocess_sol}).

\bco[Strong convergence of $Y_t$ to $X_t$]\label{NCS-corollary:convergence}
Let the assumption of Theorem \ref{NCS-theorem:convergence} hold. The semi-discrete scheme described by  (\ref{NCS-eq:SDmethod_t}) converges strongly in the mean-square sense to the true solution of (\ref{NCS-eq:systemSDEs}) with order of convergence $1/4$, that is 
 
$$\bfE\sup_{0\leq t\leq T} \| Y_t - X_t\|_2^2\leq C\D^{1/2}.$$

\eco




\section{Proofs}\label{NCS:sec:proofs}

\subsection{Proof of Theorem \ref{NCS-theorem:convergence}}\label{NCS-proof_theorem:convergence}

Denote $\bbE_t^{i+1,i}:=\wt{Y}_{t}^{i+1,i} - \wt{X}_{t}^{i+1,i}.$ Our goal is to bound $\bfE\sup_{0\leq t\leq T}(\bbE_t^{i+1,i})^2.$ We begin with moment bounds for $|\wt{X}_{t}^{i+1,i}|^p$ and $|\wt{X}_{t}^{i+1,i}|^p$ and later we estimate the local error of the proposed semi-discrete method.

\ble[Moment bounds]\label{NCS-lem:SDuniformMomentBound}
It holds that
$$
\bfE\sup_{0\leq t\leq T}\left(|\wt{Y}_{t}^{i+1,i}|^p\vee |\wt{X}_{t}^{i+1,i}|^p\right)\leq A,
$$
for any $p\in\bbR$ where $A$ is a constant.
\ele
\bpf[Proof of Lemma \ref{NCS-lem:SDuniformMomentBound}]

Set the stopping time $\tau_R:=\inf\{t\in[0,T]: \wt{Y}_{t}^{i+1,i}>R\},$ for $R>0$ with the convention $\inf \emptyset=\infty.$ Application of It\^o's formula on $(\wt{Y}_{t\wedge\tau_R}^{i+1,i})^p$, see (\ref{NCS-eq:SDprocess}) implies
\beao
(\wt{Y}_{t\wedge\tau_R}^{i+1,i})^p &=& (\wt{Y}_{0}^{i+1,i})^p +  \int_0^{t\wedge\tau_R} \left(p\alpha^i(\wt{Y}_{\hat{s}})(\wt{Y}_{s}^{i+1,i})^p - p\beta^i(\wt{Y}_{\hat{s}})(\wt{Y}_{s}^{i+1,i})^{p-1} + \frac{p(p-1)}{2} (c^i)^2 (\wt{Y}_{s}^{i+1,i})^{p-2} \right)ds\\
&&  + pc^i \int_0^{t\wedge\tau_R}  (\wt{Y}_{s}^{i+1,i})^{p-1}dB_s^{(i+1)}.
\eeao

Let us denote $b^i_{Lip}$ the Lispchitz constant of $b^i$ and $b_{Lip}=\sup_{1\leq i\leq d} b^i_{Lip}.$ 
\beam\nonumber
b^{i+1}(v^{(i+1)}) - b^{i}(v^{(i)}) &=& b^{i+1}(v^{(i+1)}) - b^{i+1}(v^{(i)})
 + b^{i+1}(v^{(i)}) - b^{i}(v^{(i)})\\
\nonumber&\leq& b^{i+1}_{Lip}v^{i+1,i} + (b^{i+1}(0) + b^{i}(0))+ (b^{i+1}_{Lip}+b^{i}_{Lip})|v^{(i)}|\\
\label{NCS-eq:b_coef0}&\leq& b_{Lip}v^{i+1,i} + (b^{i+1}(0) + b^{i}(0))+ 2b_{Lip}|v^{(i)}|,
\eeam
which implies
\beqq \label{NCS-eq:b_coef1}
\left|\frac{b^{i+1}(v^{(i+1)}) - b^{i}(v^{(i)})}{v^{i+1,i}}\right|
\leq b_{Lip} + 3\frac{|b^{i+1}(0) + b^{i}(0)|}{|v^{i+1,i}|}+ 2b_{Lip}\frac{|v^{(i)}|}{|v^{i+1,i}|}.
\eeqq

Moreover
$$\sum_{k\neq i,i+1}\frac{1}{(\wt{Y}_{\hat{s}}^{i,k})(\wt{Y}_{\hat{s}}^{i+1,k})} \leq Q^2\sup_{1\leq i\leq d-1}(\wt{Y}_{\hat{s}}^{i+1,i})^{-2},
$$
where 
\beqq \label{NCS-eq:Q_coef}
Q:=1 + 2\sum_{k=2}^i\frac{1}{k}  + 2\sum_{k=2}^{d-i}\frac{1}{k}.
\eeqq

Recall the definitions (\ref{NCS-eq:SDprocess_a}),  (\ref{NCS-eq:SDprocess_b}) of the functions $\alpha^i$ and $\beta^i.$ 
$$
|\alpha^i (\wt{Y}_{\hat{s}})| \leq  b_{Lip} + 3\frac{|b^{i+1}(0) + b^{i}(0)|}{|\wt{Y}_{\hat{s}}^{i+1,i}|}+ 2b_{Lip}\frac{|\wt{Y}_{\hat{s}}^{(i)}|}{|\wt{Y}_{\hat{s}}^{i+1,i}|} +  \gamma Q^2\sup_{1\leq i\leq d-1}(\wt{Y}_{\hat{s}}^{i+1,i})^{-2},
$$
where $\gamma := \sup_{i\neq j}\gamma_{i,j},$ and 
$$
|\beta^i (\wt{Y}_{\hat{s}})| \leq  2\gamma Q\sup_{1\leq i\leq d-1}(\wt{Y}_{\hat{s}}^{i+1,i})^{-1}.$$

Using  the inequality $x^{p-k}y\leq \ep^k\frac{p-k}{p}x^p + \frac{k}{p\ep^{p-k}}y^{p/k},$ valid for $x\wedge y\geq0$ and $p>k$ with $\ep=\frac{1}{2},$  and $k=1,2$ and the fact that $\wt{Y}_{\hat{s}}^{(i)}=\sum_{j=0}^{i-1}\wt{Y}_{\hat{s}}^{j+1,j}$
with the convection $\wt{Y}_{\hat{s}}^{1,0}\equiv\wt{Y}_{\hat{s}}^{(1)}$ we reach

$$
(\wt{Y}_{t\wedge\tau_R}^{i+1,i})^p \leq (\wt{Y}_{0}^{i+1,i})^p +  \int_0^{t\wedge\tau_R} \left(C_1 + C_2(\wt{Y}_{s}^{i+1,i})^p +C_3\left(\sum_{j=0}^{i-1}\wt{Y}_{\hat{s}}^{j+1,j}\right)^p\right)ds  + M_t,$$
where 
$$
M_t:=pc^i \int_0^{t\wedge\tau_R}  (\wt{Y}_{s}^{i+1,i})^{p-1}dB_s^{(i+1)},
$$

$$
C_1:= pb_{Lip} +\frac{3(p-1)|b^{i+1}(0) + b^{i}(0)|}{2} + (p-1)b_{Lip} + \left(\frac{(p-1)}{2} (c^i)^2 + \gamma (Q^2 + 2Q)\right)\frac{p-2}{4}
$$

$$
C_2:= 3|b^{i+1}(0) + b^{i}(0)|2^{p-1} + \left(\frac{(p-1)}{2} (c^i)^2 + \gamma (Q^2 + 2Q)\right)2^{p-1},  \quad C_3:=2^pb_{Lip}.
$$

Furthermore
\beao
\left(\sum_{j=0}^{i-1}\wt{Y}_{\hat{s}}^{j+1,j}\right)^p&\leq& i^{p-1}\sum_{j=0}^{i-1}\left(\wt{Y}_{\hat{s}}^{j+1,j}\right)^p\\
&\leq& i^{p-1}\frac{i(i+1)}{2}\sup_{0\leq i \leq d-1}\left(\wt{Y}_{s}^{i+1,i}\right)^p\\
&\leq& d^{p-1}\frac{d(d+1)}{2}\sup_{0\leq i \leq d-1}\left(\wt{Y}_{s}^{i+1,i}\right)^p.
\eeao

Thus
$$
\sup_{0\leq i \leq d-1}(\wt{Y}_{t\wedge\tau_R}^{i+1,i})^p \leq \sup_{0\leq i \leq d-1}(\wt{Y}_{0}^{i+1,i})^p +  \int_0^{t\wedge\tau_R} \left(C_1 + C_4\sup_{0\leq i \leq d-1}(\wt{Y}_{s}^{i+1,i})^p \right)ds + M_t
$$
where 
 $$
 C_4:= C_2 + C_3d^{p-1}\frac{d(d+1)}{2}.
 $$

Taking expectations in the above inequality and using that $M_t$ is a local martingale vanishing at $0,$ we get

\beao
\bfE\sup_{0\leq i \leq d-1}(\wt{Y}_{t\wedge\tau_R}^{i+1,i})^p &\leq& \bfE\sup_{0\leq i \leq d-1}(\wt{Y}_{0}^{i+1,i})^p +  \int_0^{t\wedge\tau_R} \left(C_1 + C_4\bfE\sup_{0\leq i \leq d-1}(\wt{Y}_{s}^{i+1,i})^p\right)ds \\
&\leq& \left(\bfE\sup_{0\leq i \leq d-1}(\wt{Y}_{0}^{i+1,i})^p + C_1T\right)e^{C_4T},
\eeao
where we have applied the Gronwall inequality. 
Taking the limit as $\Rto$ and applying the monotone convergence theorem leads to 
$$
\bfE\sup_{0\leq i \leq d-1}(\wt{Y}_{t}^{i+1,i})^p \leq \left(\bfE\sup_{0\leq i \leq d-1}(\wt{Y}_{0}^{i+1,i})^p + C_1T\right)e^{C_4T}.
$$
Using again It\^o's formula on $(\wt{Y}_{t}^{i+1,i})^p$, taking the supremum and then using Doob's martingale inequality on the diffusion term we bound $\bfE\sup_{0\leq t\leq T}|\wt{Y}_{t}^{i+1,i}|^p.$ The same techniques may be applied to show the result for negative $p$; the moment bounds for  $|\wt{X}_{t}^{i+1,i}|^p$ follow by similar arguments.
\epf

\ble[Local error of SD method]\label{NCS-lem:local_SD_error_p}
Let $s\in [t_{n_s}, t_{n_{s+1}}]$ where $n_s$ is an integer. Then
\beqq \label{NCS-eq:local_SD_error_p}
\bfE\sup_{0\leq i\leq d-1}|\wt{Y}_{s}^{i+1,i}-\wt{Y}_{\hat{s}}^{i+1,i}|^p\leq C\D^{p/2},
\eeqq
for any $p>0$ where the constant $C$ does not depend on $\D.$
\ele

\bpf[Proof of Lemma \ref{NCS-lem:local_SD_error_p}]\label{NCS-proof_lem:local_SD_error_p}

Take a $p>2.$ Relation (\ref{NCS-eq:SDprocess}) implies
\beao
&&|\wt{Y}_{s}^{i+1,i}-\wt{Y}_{\hat{s}}^{i+1,i}|^p=\Big|\int_{t_{n_s}}^s\left(\alpha^i(\wt{Y}_{t_{n_s}})\wt{Y}_{u}^{i+1,i} - \beta^i(\wt{Y}_{t_{n_s}}) \right)du
+ c^i \int_{t_{n_s}}^{s}  dB_s^{(i+1)}\Big|^p\\
&\leq&2^{p-1}\Big(\Big|\int_{t_{n_s}}^s\left(\alpha^i(\wt{Y}_{t_{n_s}})\wt{Y}_{u}^{i+1,i} - \beta^i(\wt{Y}_{t_{n_s}}) \right)du\Big|^p + (c^i)^p\big|\int_{t_{n_s}}^{s} dB_s^{(i+1)}\big|^p\Big)\\
&\leq&2^{p-1}\Big(|s-t_{n_s}|^{p-1}\int_{t_{n_s}}^s\left|\alpha^i(\wt{Y}_{t_{n_s}})\wt{Y}_{u}^{i+1,i} - \beta^i(\wt{Y}_{t_{n_s}}) \right|^p du + (c^i)^p\big|\int_{t_{n_s}}^{s} dB_s^{(i+1)}\big|^p\Big)\\
&\leq&4^{p-1}\D^{p-1}|\alpha^i(\wt{Y}_{t_{n_s}})|^p \int_{t_{n_s}}^s|  \wt{Y}_{u}^{i+1,i}|^p du +  4^{p-1}\D^p|\beta^i(\wt{Y}_{t_{n_s}})|^p + 2^{p-1}(c^i)^p\big|\int_{t_{n_s}}^{s} dB_s^{(i+1)}\big|^p,
\eeao
where we have used the Cauchy-Schwarz inequality. Taking expectations in the above inequality and using the uniform moment bounds of  $|\wt{Y}_{t}^{i+1,i}|^p$ described in Lemma \ref{NCS-lem:SDuniformMomentBound} and Doob's martingale inequality on the diffusion term we conclude (\ref{NCS-eq:local_SD_error_p}). The case  $0<p<2$ follows by Jensen's inequality for the concave function $\phi(x)=x^{p/2}$ since for a random variable Z it holds $\bfE|Z|^p\leq (\bfE|Z|^2)^{p/2}.$
\epf

Relations (\ref{NCS-eq:system_diff_split1}) and
(\ref{NCS-eq:SDprocess}) imply
$$
 \bbE_t^{i+1,i} = \int_0^t \left(\alpha^i(\wt{Y}_{\hat{s}})\wt{Y}_{s}^{i+1,i} - \beta^i(\wt{Y}_{\hat{s}}) +  \sum_{k\neq i,i+1}\frac{\gamma_{i,k}\wt{X}_s^{i+1,k}-\gamma_{i+1,k}\wt{X}_s^{i,k}}{\wt{X}_s^{i+1,k}\wt{X}_s^{i,k}}  - b^{i+1}(\wt{X}_s^{(i+1)}) + b^{i}(\wt{X}_s^{(i)}) \right)ds,
$$ 
where we used $\wt{Y}_{0}^{i+1,i}=\wt{X}_{0}^{i+1,i}.$ We decompose the above integrand, $I$, in the following way,
\beao
I&=& \alpha^i(\wt{Y}_{\hat{s}})(\wt{Y}_{s}^{i+1,i} - \wt{X}_{s}^{i+1,i}) + (\alpha^i(\wt{Y}_{\hat{s}}) - \alpha^i(\wt{X}_{s}) ) \wt{X}_{s}^{i+1,i} + \alpha^i(\wt{X}_s) \wt{X}_{s}^{i+1,i} - b^{i+1}(\wt{X}_s^{(i+1)}) + b^{i}(\wt{X}_s^{(i)})\\
&&+\sum_{k\neq i,i+1}\frac{\gamma_{i,k}\wt{X}_s^{i+1,k}-\gamma_{i+1,k}\wt{X}_s^{i,k}}{\wt{X}_s^{i+1,k}\wt{X}_s^{i,k}} - \beta^i(\wt{X}_{s})
+ \beta^i(\wt{X}_{s}) - \beta^i(\wt{Y}_{\hat{s}})\\
&=& \alpha^i(\wt{Y}_{\hat{s}})(\bbE_{s}^{i+1,i}) + \left(\frac{b^{i+1}(\wt{Y}_{\hat{s}}^{(i+1)})-b^{i}(\wt{Y}_{\hat{s}}^{(i)})}{\wt{Y}_{\hat{s}}^{i+1,i}} - \frac{b^{i+1}(\wt{X}_{s}^{(i+1)})-b^{i}(\wt{X}_{s}^{(i)})}{\wt{X}_{s}^{i+1,i}}\right)\wt{X}_{s}^{i+1,i}\\
&& + \left(\sum_{k\neq i,i+1}\frac{\gamma_{i,k}}{\wt{X}_{s}^{i+1,k}\wt{X}_{s}^{i,k}} - \sum_{k\neq i,i+1}\frac{\gamma_{i,k}}{\wt{Y}_{\hat{s}}^{i+1,k}\wt{Y}_{\hat{s}}^{i,k}}\right)\wt{X}_{s}^{i+1,i} - \sum_{k\neq i,i+1}\frac{\gamma_{i,k}\wt{X}_{s}^{i+1,i}}{\wt{X}_{s}^{i+1,k}\wt{X}_{s}^{i,k}}\\
&&+\sum_{k\neq i,i+1}\frac{\gamma_{i,k}\wt{X}_s^{i+1,k}-\gamma_{i+1,k}\wt{X}_s^{i,k}}{\wt{X}_s^{i+1,k}\wt{X}_s^{i,k}} -\sum_{k\neq i,i+1}\frac{\gamma_{i,k} - \gamma_{i+1,k}}{\wt{X}_s^{i+1,k}} + \sum_{k\neq i,i+1}\frac{\gamma_{i,k} - \gamma_{i+1,k}}{\wt{X}_s^{i+1,k}} - \sum_{k\neq i,i+1}\frac{\gamma_{i,k} - \gamma_{i+1,k}}{\wt{Y}_{\hat{s}}^{i+1,k}}
\eeao
where we used (\ref{NCS-eq:SDprocess_a}) and (\ref{NCS-eq:SDprocess_b}). Moreover,
\beao
&&\left(\frac{b^{i+1}(\wt{Y}_{\hat{s}}^{(i+1)})-b^{i}(\wt{Y}_{\hat{s}}^{(i)})}{\wt{Y}_{\hat{s}}^{i+1,i}} - \frac{b^{i+1}(\wt{X}_{s}^{(i+1)})-b^{i}(\wt{X}_{s}^{(i)})}{\wt{X}_{s}^{i+1,i}}\right)\wt{X}_{s}^{i+1,i}\\
&=& \left( b^{i+1}(\wt{Y}_{\hat{s}}^{(i+1)}) - b^{i}(\wt{Y}_{\hat{s}}^{(i)}) \right)   \left(\frac{\wt{X}_{s}^{i+1,i}}{\wt{Y}_{\hat{s}}^{i+1,i}}-1 \right) + b^{i+1}(\wt{Y}_{\hat{s}}^{(i+1)}) - b^{i+1}(\wt{X}_{s}^{(i+1)})  - b^{i}(\wt{Y}_{\hat{s}}^{(i)}) + b^{i}(\wt{X}_{s}^{(i)})\\
&=& \frac{ b^{i+1}(\wt{Y}_{\hat{s}}^{(i+1)}) - b^{i}(\wt{Y}_{\hat{s}}^{(i)})}{\wt{Y}_{\hat{s}}^{i+1,i}} (\wt{Y}_{s}^{i+1,i}-\wt{Y}_{\hat{s}}^{i+1,i} - \bbE_{s}^{i+1,i}) + b^{i+1}(\wt{Y}_{\hat{s}}^{(i+1)}) - b^{i+1}(\wt{X}_{s}^{(i+1)})  - b^{i}(\wt{Y}_{\hat{s}}^{(i)}) + b^{i}(\wt{X}_{s}^{(i)})
\eeao

and

\beao
&&\sum_{k\neq i,i+1}\frac{\gamma_{i,k}\wt{X}_s^{i+1,k}-\gamma_{i+1,k}\wt{X}_s^{i,k}}{\wt{X}_s^{i+1,k}\wt{X}_s^{i,k}} -\sum_{k\neq i,i+1}\frac{\gamma_{i,k} - \gamma_{i+1,k}}{\wt{X}_s^{i+1,k}}- \sum_{k\neq i,i+1}\frac{\gamma_{i,k}\wt{X}_{s}^{i+1,i}}{\wt{X}_{s}^{i+1,k}\wt{X}_{s}^{i,k}}\\
&=&\sum_{k\neq i,i+1}\left(\frac{\gamma_{i,k}}{\wt{X}_s^{i,k}}-\frac{\gamma_{i+1,k}}{\wt{X}_s^{i+1,k}} -\frac{\gamma_{i,k}}{\wt{X}_s^{i+1,k}} + \frac{\gamma_{i+1,k}}{\wt{X}_s^{i+1,k}}\right) - \sum_{k\neq i,i+1}\frac{\gamma_{i,k}\wt{X}_{s}^{i+1,i}}{\wt{X}_{s}^{i+1,k}\wt{X}_{s}^{i,k}}\\
&=&\sum_{k\neq i,i+1}\gamma_{i,k}\left(\frac{1}{\wt{X}_s^{i,k}} -\frac{1}{\wt{X}_s^{i+1,k}} - \frac{\wt{X}_{s}^{i+1,i}}{\wt{X}_{s}^{i+1,k}\wt{X}_{s}^{i,k}}\right)=0,
\eeao
since $\wt{X}_{s}^{i+1,k}-\wt{X}_{s}^{i,k}=\wt{X}_{s}^{i+1,i}$ and 

\beao
&& \left(\sum_{k\neq i,i+1}\frac{\gamma_{i,k}}{\wt{X}_{s}^{i+1,k}\wt{X}_{s}^{i,k}} - \sum_{k\neq i,i+1}\frac{\gamma_{i,k}}{\wt{Y}_{\hat{s}}^{i+1,k}\wt{Y}_{\hat{s}}^{i,k}}\right)\wt{X}_{s}^{i+1,i} + \sum_{k\neq i,i+1}\frac{\gamma_{i,k} - \gamma_{i+1,k}}{\wt{X}_s^{i+1,k}} - \sum_{k\neq i,i+1}\frac{\gamma_{i,k} - \gamma_{i+1,k}}{\wt{Y}_{\hat{s}}^{i+1,k}}\\
&=&\sum_{k\neq i,i+1}\frac{\gamma_{i,k}\wt{X}_{s}^{i+1,i}}{\wt{X}_{s}^{i+1,k}\wt{X}_{s}^{i,k}} - \sum_{k\neq i,i+1}\frac{\gamma_{i,k}\wt{X}_{s}^{i+1,i}}{\wt{Y}_{\hat{s}}^{i+1,k}\wt{Y}_{\hat{s}}^{i,k}} + \sum_{k\neq i,i+1}(\gamma_{i,k} - \gamma_{i+1,k})\left(\frac{1}{\wt{X}_s^{i+1,k}} -\frac{1}{\wt{Y}_{\hat{s}}^{i+1,k}}\right)\\
&=&\sum_{k\neq i,i+1}\gamma_{i,k}\left(\frac{1}{\wt{X}_{s}^{i,k}} -\frac{1}{\wt{Y}_{\hat{s}}^{i,k}}\right) + \sum_{k\neq i,i+1}\gamma_{i,k}\left(\frac{1}{\wt{Y}_{\hat{s}}^{i,k}} -\frac{1}{\wt{Y}_{\hat{s}}^{i+1,k}} - \frac{\wt{X}_{s}^{i+1,i}}{\wt{Y}_{\hat{s}}^{i+1,k}\wt{Y}_{\hat{s}}^{i,k}}\right)\\
&&- \sum_{k\neq i,i+1}\gamma_{i+1,k}\left(\frac{1}{\wt{X}_{s}^{i+1,k}} -\frac{1}{\wt{Y}_{\hat{s}}^{i+1,k}}\right)\\
&=&\sum_{k\neq i,i+1}\gamma_{i,k}\left(\frac{1}{\wt{X}_{s}^{i,k}} -\frac{1}{\wt{Y}_{\hat{s}}^{i,k}}\right) - \gamma_{i+1,k}\left(\frac{1}{\wt{X}_{s}^{i+1,k}} -\frac{1}{\wt{Y}_{\hat{s}}^{i+1,k}}\right) + \sum_{k\neq i,i+1}\gamma_{i,k}\frac{\wt{Y}_{\hat{s}}^{i+1,i}-\wt{X}_{s}^{i+1,i}}{\wt{Y}_{\hat{s}}^{i+1,k}\wt{Y}_{\hat{s}}^{i,k}}\\
&=&\sum_{k\neq i,i+1}\gamma_{i,k}\left(\frac{1}{\wt{X}_{s}^{i,k}} -\frac{1}{\wt{Y}_{\hat{s}}^{i,k}}\right) - \gamma_{i+1,k}\left(\frac{1}{\wt{X}_{s}^{i+1,k}} -\frac{1}{\wt{Y}_{\hat{s}}^{i+1,k}}\right)\\
&& + \frac{\bbE_{s}^{i+1,i} + (\wt{Y}_{\hat{s}}^{i+1,i}-\wt{Y}_{s}^{i+1,i})}{\wt{Y}_{\hat{s}}^{i+1,i}}\sum_{k\neq i,i+1}\gamma_{i,k}\frac{\wt{Y}_{\hat{s}}^{i+1,k}-\wt{Y}_{\hat{s}}^{i,k}}{\wt{Y}_{\hat{s}}^{i+1,k}\wt{Y}_{\hat{s}}^{i,k}}.
\eeao
The integrand $I$ becomes
\beao
I &=& \alpha^i(\wt{Y}_{\hat{s}})(\bbE_{s}^{i+1,i}) + \frac{ b^{i+1}(\wt{Y}_{\hat{s}}^{(i+1)}) - b^{i}(\wt{Y}_{\hat{s}}^{(i)})}{\wt{Y}_{\hat{s}}^{i+1,i}} (\wt{Y}_{s}^{i+1,i}-\wt{Y}_{\hat{s}}^{i+1,i} - \bbE_{s}^{i+1,i})\\
&& + b^{i+1}(\wt{Y}_{\hat{s}}^{(i+1)}) - b^{i+1}(\wt{X}_{s}^{(i+1)})  - b^{i}(\wt{Y}_{\hat{s}}^{(i)}) + b^{i}(\wt{X}_{s}^{(i)})\\
&& + \sum_{k\neq i,i+1}\gamma_{i,k}\left(\frac{1}{\wt{X}_{s}^{i,k}} -\frac{1}{\wt{Y}_{\hat{s}}^{i,k}}\right) - \gamma_{i+1,k}\left(\frac{1}{\wt{X}_{s}^{i+1,k}} -\frac{1}{\wt{Y}_{\hat{s}}^{i+1,k}}\right)\\
&& + \frac{\bbE_{s}^{i+1,i} + (\wt{Y}_{\hat{s}}^{i+1,i}-\wt{Y}_{s}^{i+1,i})}{\wt{Y}_{\hat{s}}^{i+1,i}}\sum_{k\neq i,i+1}\gamma_{i,k}\frac{\wt{Y}_{\hat{s}}^{i+1,k}-\wt{Y}_{\hat{s}}^{i,k}}{\wt{Y}_{\hat{s}}^{i+1,k}\wt{Y}_{\hat{s}}^{i,k}}\\
&=&  \frac{ b^{i+1}(\wt{Y}_{\hat{s}}^{(i+1)}) - b^{i}(\wt{Y}_{\hat{s}}^{(i)})}{\wt{Y}_{\hat{s}}^{i+1,i}} (\wt{Y}_{s}^{i+1,i}-\wt{Y}_{\hat{s}}^{i+1,i})\\
&& + b^{i+1}(\wt{Y}_{\hat{s}}^{(i+1)}) - b^{i+1}(\wt{X}_{s}^{(i+1)})  - b^{i}(\wt{Y}_{\hat{s}}^{(i)}) + b^{i}(\wt{X}_{s}^{(i)})\\
&& + \sum_{k\neq i,i+1}\gamma_{i,k}\left(\frac{1}{\wt{X}_{s}^{i,k}} -\frac{1}{\wt{Y}_{\hat{s}}^{i,k}}\right) - \gamma_{i+1,k}\left(\frac{1}{\wt{X}_{s}^{i+1,k}} -\frac{1}{\wt{Y}_{\hat{s}}^{i+1,k}}\right)\\
&& + (\wt{Y}_{\hat{s}}^{i+1,i}-\wt{Y}_{s}^{i+1,i})\sum_{k\neq i,i+1}\frac{\gamma_{i,k}}{\wt{Y}_{\hat{s}}^{i+1,k}\wt{Y}_{\hat{s}}^{i,k}}
\eeao
where we used once more (\ref{NCS-eq:SDprocess_a}).
\beam\nonumber
(\bbE_t^{i+1,i})^2 &\leq&5t\int_0^t \left(\frac{ b^{i+1}(\wt{Y}_{\hat{s}}^{(i+1)}) - b^{i}(\wt{Y}_{\hat{s}}^{(i)})}{\wt{Y}_{\hat{s}}^{i+1,i}}\right)^2 (\wt{Y}_{s}^{i+1,i}-\wt{Y}_{\hat{s}}^{i+1,i})^2ds\\
\nonumber&&+ 5t\int_0^t  \left(b^{i+1}(\wt{Y}_{\hat{s}}^{(i+1)}) - b^{i+1}(\wt{X}_{s}^{(i+1)})  - b^{i}(\wt{Y}_{\hat{s}}^{(i)}) + b^{i}(\wt{X}_{s}^{(i)})\right)^2ds\\
\nonumber&&+ 5t\int_0^t \left(\sum_{k\neq i.i+1}\gamma_{i,k}\left(\frac{1}{\wt{X}_{s}^{i,k}} -\frac{1}{\wt{Y}_{\hat{s}}^{i,k}}\right)\right)^2 +\left(\sum_{k\neq i,i+1}\gamma_{i+1,k}\left(\frac{1}{\wt{X}_{s}^{i+1,k}} -\frac{1}{\wt{Y}_{\hat{s}}^{i+1,k}}\right)\right)^2ds\\
\label{NCS-eq:dist1}&&+  5t\int_0^t (\wt{Y}_{s}^{i+1,i}-\wt{Y}_{\hat{s}}^{i+1,i})^2\left(\sum_{k\neq i,i+1}\frac{\gamma_{i,k}}{\wt{Y}_{\hat{s}}^{i,k}\wt{Y}_{\hat{s}}^{i+1,k}}\right)^2ds. 
\eeam

Squaring both sides of inequality (\ref{NCS-eq:b_coef0}) yields 
\beqq \label{NCS-eq:b_coef}
\left(\frac{b^{i+1}(v^{(i+1)}) - b^{i}(v^{(i)})}{v^{i+1,i}}\right)^2
\leq 3b_{Lip}^2 + 3\frac{(b^{i+1}(0) + b^{i}(0))^2}{(v^{i+1,i})^2}+ 4b^2_{Lip}\frac{(v^{(i)})^2}{(v^{i+1,i})^2}.
\eeqq

Furthermore,
\beao
&&\left(\sum_{k\neq i,i+1}\gamma_{i,k}\left(\frac{1}{\wt{X}_{s}^{i,k}} -\frac{1}{\wt{Y}_{\hat{s}}^{i,k}}\right)\right)^2 \leq  2\left(\sum_{k\neq i,i+1}\gamma_{i,k}\frac{\wt{Y}_{s}^{i,k}-\wt{Y}_{\hat{s}}^{i,k}}{\wt{X}_{s}^{i,k}\wt{Y}_{\hat{s}}^{i,k}}\right)^2
+ 2\left(\sum_{k\neq i,i+1}\gamma_{i,k}\frac{\bbE_{s}^{i,k}}{\wt{X}_{s}^{i,k}\wt{Y}_{\hat{s}}^{i,k}}\right)^2\\
&\leq&  2\hat{Q}^2\sup_{0\leq i\leq d-1} (\wt{Y}_{s}^{i+1,i}-\wt{Y}_{\hat{s}}^{i+1,i})^2
\left(\sum_{k\neq i,i+1}\frac{\gamma_{i,k}}{\wt{X}_{s}^{i,k}\wt{Y}_{\hat{s}}^{i,k}}\right)^2
+ 2\hat{Q}^2\sup_{0\leq i\leq d-1} (\bbE_{s}^{i+1,i})^2\left(\sum_{k\neq i,i+1}\frac{\gamma_{i,k}}{\wt{X}_{s}^{i,k}\wt{Y}_{\hat{s}}^{i,k}}\right)^2
\eeao
and as a consequence 
\beao
\left(\sum_{k\neq i,i+1}\gamma_{i+1,k}\left(\frac{1}{\wt{X}_{s}^{i+1,k}} -\frac{1}{\wt{Y}_{\hat{s}}^{i+1,k}}\right)\right)^2  &\leq&  \hat{Q}^2\sup_{0\leq i\leq d-1} (\wt{Y}_{s}^{i+1,i}-\wt{Y}_{\hat{s}}^{i+1,i})^2\left(\sum_{k\neq i,i+1}\frac{\gamma_{i+1,k}}{\wt{X}_{s}^{i+1,k}\wt{Y}_{\hat{s}}^{i+1,k}}\right)^2\\  
&& + \hat{Q}^2\sup_{0\leq i\leq d-1} (\bbE_{s}^{i+1,i})^2\left(\sum_{k\neq i,i+1}\frac{\gamma_{i+1,k}}{\wt{X}_{s}^{i+1,k}\wt{Y}_{\hat{s}}^{i+1,k}}\right)^2,
 \eeao
where
\beqq \label{NCS-eq:Qhat_coef}
\hat{Q}:=1 + 2\sum_{k=2}^i k   + 2\sum_{k=2}^{d-i}k.
\eeqq
Finally
\beao
&&\left(b^{i+1}(\wt{Y}_{\hat{s}}^{(i+1)}) - b^{i+1}(\wt{X}_{s}^{(i+1)})  - b^{i}(\wt{Y}_{\hat{s}}^{(i)}) + b^{i}(\wt{X}_{s}^{(i)})\right)^2\\
&\leq& 2\left(b^{i+1}(\wt{Y}_{\hat{s}}^{(i+1)}) - b^{i+1}(\wt{X}_{s}^{(i+1)})\right)^2 + 2\left( b^{i}(\wt{Y}_{\hat{s}}^{(i)}) - b^{i}(\wt{X}_{s}^{(i)})\right)^2\\
&\leq& 2(b_{Lip}^{i+1})^2(\wt{Y}_{\hat{s}}^{(i+1)} -\wt{X}_{s}^{(i+1)})^2 + 2(b_{Lip}^{i})^2(\wt{Y}_{\hat{s}}^{(i)} - \wt{X}_{s}^{(i)})^2\\
&\leq& 2b_{Lip}^2(\bbE_{s}^{(i+1)})^2 + 2b_{Lip}^2(\bbE_{s}^{(i)})^2\\
&\leq& 4b_{Lip}^2(\bbE_{s}^{i+1,i})^2 + 6b_{Lip}^2\sum_{j=0}^{i-1}(\bbE_{s}^{j+1,j})^2 \leq  4b_{Lip}^2(\bbE_{s}^{i+1,i})^2 + 6db_{Lip}^2\sup_{0\leq i\leq d-1}(\bbE_{s}^{i+1,i})^2.
\eeao

Plugging all the above estimates and (\ref{NCS-eq:b_coef}) into (\ref{NCS-eq:dist1}) yields

\beam\nonumber
(\bbE_t^{i+1,i})^2 &\leq&5t\int_0^t \left(3b_{Lip}^2 + 3\frac{(b^{i+1}(0) + b^{i}(0))^2}{(\wt{Y}_{\hat{s}}^{i+1,i})^2}+ 4b^2_{Lip}\frac{(\wt{Y}_{\hat{s}}^{(i)})^2}{(\wt{Y}_{\hat{s}}^{i+1,i})^2}\right) (\wt{Y}_{s}^{i+1,i}-\wt{Y}_{\hat{s}}^{i+1,i})^2ds\\
\nonumber&&+ 5t\int_0^t  \left(4b_{Lip}^2(\bbE_{s}^{i+1,i})^2 + 6db_{Lip}^2\sup_{0\leq i\leq d-1}(\bbE_{s}^{i+1,i})^2\right)ds\\
\nonumber&&+ 10t\int_0^t \hat{Q}^2\sup_{0\leq i\leq d-1} (\bbE_{s}^{i+1,i})^2\left[\left(\sum_{k\neq i,i+1}\frac{\gamma_{i,k}}{\wt{X}_{s}^{i,k}\wt{Y}_{\hat{s}}^{i,k}}\right)^2 +  \left(\sum_{k\neq i,i+1}\frac{\gamma_{i+1,k}}{\wt{X}_{s}^{i+1,k}\wt{Y}_{\hat{s}}^{i+1,k}}\right)^2\right]ds\\
\nonumber&&+ 10t\int_0^t \hat{Q}^2\sup_{0\leq i\leq d-1} (\wt{Y}_{s}^{i+1,i}-\wt{Y}_{\hat{s}}^{i+1,i})^2\left[\left(\sum_{k\neq i,i+1}\frac{\gamma_{i,k}}{\wt{X}_{s}^{i,k}\wt{Y}_{\hat{s}}^{i,k}}\right)^2 +  \left(\sum_{k\neq i,i+1}\frac{\gamma_{i+1,k}}{\wt{X}_{s}^{i+1,k}\wt{Y}_{\hat{s}}^{i+1,k}}\right)^2\right]ds\\
\label{NCS-eq:dist2}&&\quad+  10t\int_0^t (\wt{Y}_{s}^{i+1,i}-\wt{Y}_{\hat{s}}^{i+1,i})^2\sum_{k\neq i,i+1}\frac{\gamma_{i,k}^2}{(\wt{Y}_{\hat{s}}^{i,k})^2(\wt{Y}_{\hat{s}}^{i+1,k})^2}ds. 
\eeam

We define the process
$$
\zeta(t):=\int_0^t 10T\hat{Q}^2\gamma^2\left[\left(\sum_{k\neq i,i+1}\frac{1}{\wt{X}_{s}^{i,k}\wt{Y}_{\hat{s}}^{i,k}}\right)^2 +  \left(\sum_{k\neq i,i+1}\frac{1}{\wt{X}_{s}^{i+1,k}\wt{Y}_{\hat{s}}^{i+1,k}}\right)^2\right]ds
$$
and the stopping time
$$
\tau_l:=\inf\{s\in[0,T]: (20  + 30d)Tb_{Lip}^2  s + \zeta(s)\geq l\}.
$$

Relation (\ref{NCS-eq:dist2}) becomes

\beam\nonumber
&&\bfE\sup_{0\leq t\leq \tau}\sup_{0\leq i\leq d-1}(\bbE_t^{i+1,i})^2 \leq 15Tb_{Lip}^2\int_0^\tau \bfE\sup_{0\leq i\leq d-1}(\wt{Y}_{s}^{i+1,i}-\wt{Y}_{\hat{s}}^{i+1,i})^2ds\\
\nonumber
&&+3(b^{i+1}(0) + b^{i}(0))^2\int_0^\tau \bfE\sup_{0\leq i\leq d-1} (\wt{Y}_{\hat{s}}^{i+1,i})^{-2}(\wt{Y}_{s}^{i+1,i}-\wt{Y}_{\hat{s}}^{i+1,i})^2ds\\   
\nonumber
&&+4b^2_{Lip}\int_0^\tau \bfE\sup_{0\leq i\leq d-1}\frac{(\wt{Y}_{\hat{s}}^{(i)})^2}{(\wt{Y}_{\hat{s}}^{i+1,i})^2} (\wt{Y}_{s}^{i+1,i}-\wt{Y}_{\hat{s}}^{i+1,i})^2ds\\
\nonumber&&+ (20  + 30d)Tb_{Lip}^2\int_0^\tau  \bfE\sup_{0\leq l\leq s}\sup_{0\leq i\leq d-1}(\bbE_{l}^{i+1,i})^2ds + \int_0^\tau  (\zeta_s)^\prime\bfE\sup_{0\leq l\leq s}\sup_{0\leq i\leq d-1}(\bbE_{l}^{i+1,i})^2 ds\\
\nonumber&&+ 20T\hat{Q}^2\gamma^2\int_0^\tau \bfE\sup_{0\leq i\leq d-1} (\wt{Y}_{s}^{i+1,i}-\wt{Y}_{\hat{s}}^{i+1,i})^2\left(\sum_{k\neq i,i+1}\frac{1}{(\wt{X}_{s}^{i,k})^2(\wt{Y}_{\hat{s}}^{i,k})^2} +  \sum_{k\neq i,i+1}\frac{1}{(\wt{X}_{s}^{i+1,k})^2(\wt{Y}_{\hat{s}}^{i+1,k})^2}\right)ds\\
\label{NCS-eq:dist3}&&\qquad\qquad\qquad+  10T\gamma^2\int_0^\tau \bfE\sup_{0\leq i\leq d-1}(\wt{Y}_{s}^{i+1,i}-\wt{Y}_{\hat{s}}^{i+1,i})^2\sum_{k\neq i,i+1}\frac{1}{(\wt{Y}_{\hat{s}}^{i,k})^2(\wt{Y}_{\hat{s}}^{i+1,k})^2}ds, 
\eeam
where $\tau$ is a stopping time. The local error of the semi-discrete method implies 
\beqq \label{NCS-eq:local_SD_error}
\bfE\sup_{0\leq i\leq d-1}(\wt{Y}_{s}^{i+1,i}-\wt{Y}_{\hat{s}}^{i+1,i})^p\leq \D^{p/2}
\eeqq
and  
$$
\sum_{k\neq i,i+1}\frac{1}{(\wt{X}_{s}^{i,k})^2(\wt{Y}_{\hat{s}}^{i,k})^2}+\sum_{k\neq i,i+1}\frac{1}{(\wt{X}_{s}^{i+1,k})^2(\wt{Y}_{\hat{s}}^{i+1,k})^2}\leq 2Q^2\sup_{0\leq i\leq d-1} (\wt{Y}_{\hat{s}}^{i+1,i})^{-2}\sup_{0\leq i\leq d-1}(\wt{X}_{s}^{i+1,i})^{-2}$$
and 
$$\sum_{k\neq i,i+1}\frac{1}{(\wt{Y}_{\hat{s}}^{i,k})^2(\wt{Y}_{\hat{s}}^{i+1,k})^2} \leq Q^2\sup_{0\leq i\leq d-1}(\wt{Y}_{\hat{s}}^{i+1,i})^{-4}.
$$
We insert these bounds and (\ref{NCS-eq:local_SD_error}) to (\ref{NCS-eq:dist3}) and get 
\beam\nonumber
&&\bfE\sup_{0\leq t\leq \tau}\sup_{0\leq i\leq d-1}(\bbE_t^{i+1,i})^2 \leq C\D + C\int_0^\tau \sqrt{\bfE\sup_{0\leq i\leq d-1} (\wt{Y}_{\hat{s}}^{i+1,i})^{-4}}\sqrt{\bfE (\wt{Y}_{s}^{i+1,i}-\wt{Y}_{\hat{s}}^{i+1,i})^4}ds\\   
\nonumber
&&+C\int_0^\tau \left(\bfE\sup_{0\leq i\leq d-1}(\wt{Y}_{\hat{s}}^{(i)})^6\right)^{1/3}
\left(\bfE\sup_{0\leq i\leq d-1}(\wt{Y}_{\hat{s}}^{i+1,i})^{-6}\right)^{1/3}
\left(\bfE\sup_{0\leq i\leq d-1}(\wt{Y}_{s}^{i+1,i}-\wt{Y}_{\hat{s}}^{i+1,i})^6\right)^{1/3}ds\\
\nonumber&&+ (20  + 30d)Tb_{Lip}^2\int_0^\tau  \bfE\sup_{0\leq l\leq s}\sup_{0\leq i\leq d-1}(\bbE_{l}^{i+1,i})^2ds + \int_0^\tau  (\zeta_s)^\prime\bfE\sup_{0\leq l\leq s}\sup_{0\leq i\leq d-1}(\bbE_{l}^{i+1,i})^2 ds\\
\nonumber&&+ C\int_0^\tau \left(\bfE\sup_{0\leq i\leq d-1}(\wt{Y}_{\hat{s}}^{i+1,i})^{-6}\right)^{1/3}
\left(\bfE\sup_{0\leq i\leq d-1}(\wt{X}_{\hat{s}}^{i+1,i})^{-6}\right)^{1/3}
\left(\bfE\sup_{0\leq i\leq d-1}(\wt{Y}_{s}^{i+1,i}-\wt{Y}_{\hat{s}}^{i+1,i})^6\right)^{1/3}ds\\
\nonumber&&+ C\int_0^\tau  \left(\bfE\sup_{0\leq i\leq d-1}(\wt{Y}_{\hat{s}}^{i+1,i})^{-6}\right)^{2/3}\left(\bfE\sup_{0\leq i\leq d-1} (\wt{Y}_{s}^{i+1,i}-\wt{Y}_{\hat{s}}^{i+1,i})^6\right)^{1/3}ds\\
\nonumber&\leq& C\D + (20  + 30d)Tb_{Lip}^2\int_0^\tau  \bfE\sup_{0\leq l\leq s}\sup_{0\leq i\leq d-1}(\bbE_{l}^{i+1,i})^2ds + \int_0^\tau  (\zeta_s)^\prime\bfE\sup_{0\leq l\leq s}\sup_{0\leq i\leq d-1}(\bbE_{l}^{i+1,i})^2 ds\\
\label{NCS-eq:unifBound until stopping_time}&&\quad\qquad \leq C\D + \int_0^\tau  \left((20  + 30d)Tb_{Lip}^2  s + \zeta_s\right)^\prime\bfE\sup_{0\leq l\leq s}\sup_{0\leq i\leq d-1}(\bbE_{l}^{i+1,i})^2 ds.
\eeam

The uniform moment bound (\ref{NCS-eq:unifBound until stopping_time}) for $\tau=\tau_l$ reads
\beam\nonumber
\bfE\sup_{0\leq t\leq \tau_l}\sup_{0\leq i\leq d-1}(\bbE_t^{i+1,i})^2 &\leq& C\D + \int_0^{\tau_l}  \left((20  + 30d)Tb_{Lip}^2  s + \zeta_s\right)^\prime\bfE\sup_{0\leq l\leq s}\sup_{0\leq i\leq d-1}(\bbE_{l}^{i+1,i})^2 ds\\
\nonumber&\leq& C\D  + \int_0^l \bfE\sup_{0\leq j\leq u}\sup_{0\leq i\leq d-1}(\bbE_{\tau_j}^{i+1,i})^2du\\
\label{NCS-eq:unifBound until tau} &\leq& C\D e^l,
\eeam
where in the final step we have used Gronwall's inequality. Under the change of variables $u=(20  + 30d)Tb_{Lip}^2  s + \zeta_s$ relation (\ref{NCS-eq:unifBound until stopping_time})  for $\tau=T$ becomes 
\beam\nonumber
&&\bfE\sup_{0\leq t\leq T}\sup_{0\leq i\leq d-1}(\bbE_t^{i+1,i})^2 \leq C\D + \int_0^{(20  + 30d)b_{Lip}^2  T^2 + \zeta_T}\bfE\sup_{0\leq j\leq u}\sup_{0\leq i\leq d-1}(\bbE_{\tau_j}^{i+1,i})^2du \\
\nonumber&\leq&C\D + \int_0^{\infty}\bfE\left(\sup_{0\leq j\leq u}\sup_{0\leq i\leq d-1}(\bbi_{(20  + 30d)b_{Lip}^2  T^2 + \zeta_T\geq u\}}(\bbE_{\tau_j}^{i+1,i})^2\right)du \\
\nonumber&\leq&C\D + \int_0^{(20  + 30d)b_{Lip}^2  T^2}\bfE\sup_{0\leq j\leq u}\sup_{0\leq i\leq d-1}(\bbE_{\tau_j}^{i+1,i})^2du\\
\nonumber&& + \int_{(20  + 30d)b_{Lip}^2  T^2}^{\infty}\bfP((20  + 30d)b_{Lip}^2  T^2 + \zeta_T\geq u)\\
\nonumber&&\qquad\times\bfE\left(\sup_{0\leq j\leq u}\sup_{0\leq i\leq d-1}(\bbE_{\tau_j}^{i+1,i})^2 \big| \{(20  + 30d)b_{Lip}^2  T^2 + \zeta_T\geq u\}\right)du\\
\nonumber&\leq&C\D + Ce^{(20  + 30d)b_{Lip}^2  T^2}\D + \int_0^{\infty}\bfP(\zeta_T\geq u)\bfE\sup_{0\leq j\leq u}\sup_{0\leq i\leq d-1}(\bbE_{\tau_j}^{i+1,i})^2 du \\
\label{NCS-eq:unifBound until T}&\leq&C\D +C\D\int_0^{\infty}\bfP(\zeta_T\geq u)e^udu,
\eeam
where in the last steps we have used (\ref{NCS-eq:unifBound until tau}). The next step is to show $u\rightarrow \bfP(\zeta_T\geq u)e^u\in\bbl^1(\bbR_+).$ Markov's inequality implies 
$$ 
\bfP(\zeta_T\geq u)\leq e^{-\ep u}\bfE(e^{\ep \zeta_T}),
$$ 
for any $\ep>0.$ The following bound holds
\beao
\zeta_T&\leq&\int_0^T 10T\hat{Q}^2\gamma^2\left[\left(\sum_{k\neq i,i+1}\frac{1}{\wt{X}_{s}^{i,k}\wt{Y}_{\hat{s}}^{i,k}}\right)^2 +  \left(\sum_{k\neq i,i+1}\frac{1}{\wt{X}_{s}^{i+1,k}\wt{Y}_{\hat{s}}^{i+1,k}}\right)^2\right]ds\\
&\leq&20T\hat{Q}^2Q^2\gamma^2\int_0^T  \sup_{0\leq i\leq d-1} (\wt{Y}_{\hat{s}}^{i+1,i})^{-2}\sup_{0\leq i\leq d-1}(\wt{X}_{s}^{i+1,i})^{-2}ds,
\eeao
therefore
\beqq \label{NCS-eq:exponentialMoment}
\bfE(e^{\ep\zeta_T})\leq \bfE\left(e^{\ep 20T\hat{Q}^2Q^2\gamma^2\int_0^T \sup_{0\leq i\leq d-1} (\wt{Y}_{\hat{s}}^{i+1,i})^{-2}\sup_{0\leq i\leq d-1}(\wt{X}_{s}^{i+1,i})^{-2}ds} \right).
\eeqq

Applying the It\^o formula to (\ref{NCS-eq:system_diff_split1})

\beao
(\wt{X}_t^{i+1,i})^2 &=& (\wt{X}_0^{i+1,i})^2 + \int_0^t  2\wt{X}_s^{i+1,i}\left( - \sum_{k\neq i,i+1}\frac{\gamma_{i,k}\wt{X}_s^{i+1,k}-\gamma_{i+1,k}\wt{X}_s^{i,k}}{\wt{X}_s^{i+1,k}\wt{X}_s^{i,k}}\right)ds\\
&&  +\int_0^t \left[2\wt{X}_s^{i+1,i}\left(b^{i+1}(\wt{X}_s^{(i+1)})-b^{i}(\wt{X}_s^{(i)})  \right) + (c^{i})^2\right]ds + 2(c^{i})\int_0^t \wt{X}_s^{i+1,i}dB_s^{(i+1)}\\
&\geq& (\wt{X}_0^{i+1,i})^2 + \int_0^t  2\wt{X}_s^{i+1,i}\left( - \sum_{k\neq i,i+1}\frac{\gamma_{i,k}}{\wt{X}_s^{i,k}} + \sum_{k\neq i,i+1}\frac{\gamma_{i+1,k}}{\wt{X}_s^{i+1,k}}\right)ds\\
&&  +\int_0^t 2\left[\wt{X}_s^{i+1,i}\left(b^{i}(\wt{X}_s^{(i+1)})-b^{i}(\wt{X}_s^{(i)})  \right) + (c^{i})^2\right]ds + 2(c^{i})\int_0^t \wt{X}_s^{i+1,i}dB_s^{(i+1)}\\
&\geq& (\wt{X}_0^{i+1,i})^2 + \int_0^t  \left(-Q\gamma  -2b_{Lip}(\wt{X}_s^{i+1,i})^2 + (c^{i})^2\right)ds + 2(c^{i})\int_0^t \wt{X}_s^{i+1,i}dB_s^{(i+1)}.
\eeao
In the event $(\wt{X}_t^{i+1,i})^2\geq1$ we have $(\wt{X}_t^{i+1,i})^{-2}\leq1$ whereas when $(\wt{X}_t^{i+1,i})^2\leq1$ and 
\beqq \label{NCS-eq:param_cond}
(c^{i})^2 \geq Q\gamma,  
\eeqq
we get
\beao
(\wt{X}_t^{i+1,i})^2 &\geq& (\wt{X}_0^{i+1,i})^2 + \int_0^t -2b_{Lip}\wt{X}_s^{i+1,i}ds + 2(c^{i})\int_0^t \wt{X}_s^{i+1,i}dB_s^{(i+1)}\\
&\geq& (\wt{X}_0^{i+1,i})^2\exp\left\{ \int_0^t \left(-2b_{Lip}-2(c^{i})^2\right) ds + 2(c^{i})\int_0^t dB_s^{(i+1)}\right\}.
\eeao
Consequently,
\beam\nonumber
(\wt{X}_t^{i+1,i})^{-2} &\leq& (\wt{X}_0^{i+1,i})^{-2}\exp\left\{ \int_0^t \left(2b_{Lip}+2(c^{i})^2\right) ds - 2(c^{i})\int_0^t dB_s^{(i+1)}\right\}\\
\nonumber&\leq& (\wt{X}_0^{i+1,i})^{-2}e^{(2b_{Lip}+4(c^{i})^2)T}\exp\left\{ \int_0^t -2(c^{i})^2 ds - 2(c^{i})\int_0^t dB_s^{(i+1)}\right\}\\
\label{NCS-eq:X_inv_bound}&\leq& (\wt{X}_0^{i+1,i})^{-2}e^{(2b_{Lip}+4(c^{i})^2)T}\xi_t
\eeam
where $\xi_t$ is the exponential martingale

\beqq \label{NCS-eq:exp_mart}
\xi_t:=\exp\left\{ \int_0^t -2(c^{i})^2 ds - 2(c^{i})\int_0^t dB_s^{(i+1)}\right\}.\eeqq

Note that (\ref{NCS-eq:param_cond}) holds when 
\beqq \label{NCS-eq:param_cond_d}
(c^{i})^2 \geq (d-1)\gamma,  
\eeqq
since by the definition of $Q$ (\ref{NCS-eq:Q_coef}) we get
$$
Q\leq 1 + 2\frac{i-1}{2}  + 2\frac{d-i+1}{2}\leq  d-1.
$$

In the same spirit we bound $(\wt{Y}_t^{i+1,i})^{-2}.$ Once again we apply  It\^o's formula to (\ref{NCS-eq:SDprocess})

$$(\wt{Y}_{t}^{i+1,i})^2  = (\wt{Y}_{0}^{i+1,i})^2 +  \int_0^t \left(2\alpha^i(\wt{Y}_{\hat{s}})(\wt{Y}_{s}^{i+1,i})^2 - 2\wt{Y}_{s}^{i+1,i}\beta^i(\wt{Y}_{\hat{s}}) + ( c^i)^2 \right)ds  + 2c^i \int_0^t  \wt{Y}_s^{i+1,i}dB_s^{(i+1)}.
$$

Moreover,
\beao
&&2\alpha^i(\wt{Y}_{\hat{s}})(\wt{Y}_{s}^{i+1,i})^2 - 2\wt{Y}_{s}^{i+1,i}\beta^i(\wt{Y}_{\hat{s}}) \geq 2 \frac{b^{i}(\wt{Y}_{\hat{s}}^{(i+1)})-b^{i}(\wt{Y}_{\hat{s}}^{(i)})}{\wt{Y}_{s}^{i+1,i}}(\wt{Y}_{s}^{i+1,i})^2\\
&&  -\wt{Y}_{s}^{i+1,i}\sum_{k\neq i,i+1}2\gamma_{i,k}\frac{ \wt{Y}_{\hat{s}}^{i+1,k}-\wt{Y}_{\hat{s}}^{i,k}}{\wt{Y}_{\hat{s}}^{i,k}\wt{Y}_{\hat{s}}^{i+1,k}} -2\wt{Y}_s^{i+1,i}\sum_{k\neq i,i+1}\frac{\gamma_{i,k} - \gamma_{i+1,k}}{\wt{Y}_{\hat{s}}^{i+1,k}}\\
&\geq& -2 b_{Lip}(\wt{Y}_{s}^{i+1,i})^2 -2\wt{Y}_s^{i+1,i}\left(\sum_{k\neq i,i+1}\frac{\gamma_{i,k}}{\wt{Y}_{\hat{s}}^{i,k}}-\sum_{k\neq i,i+1}\frac{\gamma_{i+1,k}}{\wt{Y}_{\hat{s}}^{i+1,k}}\right)\\
&\geq& -2 b_{Lip}(\wt{Y}_{s}^{i+1,i})^2 -Q\gamma,
\eeao
therefore by a comparison theorem we have 
$$
(\wt{Y}_{t}^{i+1,i})^2\geq (\wt{Y}_{0}^{i+1,i})^2 +  \int_0^t \left[(c^i)^2-3\gamma d(d-1) - 2b_{Lip}(\wt{Y}_s^{i+1,i})^2\right]ds + 2c^i \int_0^t  \wt{Y}_s^{i+1,i}dB_s^{(i+1)}.
$$

In the event $(\wt{Y}_t^{i+1,i})^2\geq1$ we have $(\wt{Y}_t^{i+1,i})^{-2}\leq1$ whereas when $(\wt{Y}_t^{i+1,i})^2\leq1$

\beao
(\wt{Y}_{t}^{i+1,i})^2 &\geq& (\wt{Y}_{0}^{i+1,i})^2 +  \int_0^t \left[(c^i)^2-3\gamma d(d-1) - 2b_{Lip}(\wt{Y}_s^{i+1,i})\right]ds + 2c^i \int_0^t  \wt{Y}_s^{i+1,i}dB_s^{(i+1)}\\
&\geq& (\wt{Y}_{0}^{i+1,i})^2 +  \int_0^t -2b_{Lip}(\wt{Y}_s^{i+1,i})ds + 2c^i \int_0^t  \wt{Y}_s^{i+1,i}dB_s^{(i+1)}\\
&\geq& (\wt{Y}_0^{i+1,i})^2\exp\left\{ \int_0^t \left(-2b_{Lip}-2(c^{i})^2\right) ds + 2(c^{i})\int_0^t dB_s^{(i+1)}\right\},
\eeao
when condition  (\ref{NCS-eq:param_cond}) holds.

Therefore,
\beqq\label{NCS-eq:Y_inv_bound}
(\wt{Y}_t^{i+1,i})^{-2} \leq (\wt{Y}_0^{i+1,i})^{-2}e^{(2b_{Lip}+4(c^{i})^2)T}\xi_t,
\eeqq
where $\xi_t$ is the exponential martingale of (\ref{NCS-eq:exp_mart}). Now, we plug estimates (\ref{NCS-eq:X_inv_bound}) and (\ref{NCS-eq:Y_inv_bound})  into (\ref{NCS-eq:exponentialMoment}) to get

\beao
&&\bfE(e^{\ep\zeta_T}) \leq \bfE\left(e^{\ep 20T\hat{Q}^2Q^2\gamma^2\int_0^T (\wt{X}_0^{i+1,i})^{-2}e^{(2b_{Lip}+4(c^{i})^2)T}\xi_s(\wt{Y}_0^{i+1,i})^{-2}e^{(2b_{Lip}+4(c^{i})^2)T}\xi_s ds} \right)\\
&\leq &\bfE\left(e^{\ep 20T\hat{Q}^2Q^2\gamma^2e^{(2b_{Lip}+4(c^{i})^2)2T}(\wt{X}_0^{i+1,i})^{-4}\int_0^T (\xi_s)^2 ds} \right)\\
&\leq &\bfE\left(e^{\ep 20T\hat{Q}^2Q^2\gamma^2e^{(2b_{Lip}+4(c^{i})^2)2T}\int_0^T(\wt{X}_0^{i+1,i})^{-8}+ (\xi_s)^4 ds} \right)\\
&\leq &\sqrt{\bfE\left(e^{\ep 40T^2\hat{Q}^2Q^2\gamma^2e^{(2b_{Lip}+4(c^{i})^2)2T}(\wt{X}_0^{i+1,i})^{-8}} \right)}\sqrt{\bfE\left(e^{\ep 40T\hat{Q}^2Q^2\gamma^2e^{(2b_{Lip}+4(c^{i})^2)2T}\int_0^T(\xi_s)^4 ds} \right)}\\
&\leq&C\sqrt{\bfE\left(e^{\ep 40T\hat{Q}^2Q^2\gamma^2e^{(2b_{Lip}+4(c^{i})^2)2T}\int_0^T(\xi_s)^4 ds} \right)},
\eeao
where we used the exponential inverse moment bound $\bfE e^{C(\wt{X}_0^{i+1,i})^{-8}}<A_{X_0}$ for any $C>0$ where $A_{X_0}$ is a finite  constant.
By (\ref{NCS-eq:exp_mart}) we have
\beao
(\xi_t)^4&=&\exp\left\{ \int_0^t -8(c^{i})^2 ds - 8(c^{i})\int_0^t dB_s^{(i+1)}\right\}\\
&=&e^{24(c^{i})^2t}\exp\left\{ \int_0^t -32(c^{i})^2 ds - 8(c^{i})\int_0^t dB_s^{(i+1)}\right\}\\
&\leq&e^{24(c^{i})^2T}\lam_t,
\eeao
where $\lam_t$ is the solution to the  SDE

\beqq\label{NCS-eq:exp_mart_aux}
d \lam_t = -8(c^{i})dB_t^{(i+1)}.
\eeqq

Finally $\bfE\sup_{\leq t \leq T} e^{C\lam_t}<A_{\lam},$ for any $C>0$ where $A_{\lam}$ is a positive constant, thus we can always find an $\ep>1$ such that  $\bfE(e^{\ep\zeta_T})<C$ or $\bfP(\zeta_T\geq u)\leq Ce^{-\ep u}$; we return to estimate (\ref{NCS-eq:unifBound until T}) and conclude

\beqq\label{NCS-eq:unifBound until T_final}
\bfE\sup_{0\leq t\leq T}\sup_{0\leq i\leq d-1}(\bbE_t^{i+1,i})^2  \leq C\D +C\D\int_0^{\infty}\bfP(\zeta_T\geq u)e^{(1-\ep) u}du\leq C\D,
\eeqq
by choosing  $\ep>1.$

\subsection{Proof of Corollary \ref{NCS-corollary:convergence}}\label{NCS-proof_corollary:convergence}
It holds that
\beao
 \| Y_t - X_t\|_2^2 = \sum_{i=1}^d|Y_t^{(i)} - X_t^{(i)}|^2 &=&
 \sum_{i=1}^d \left| \sum_{j=0}^{i-1} Y_t^{j+1,j}  - \sum_{j=0}^d X_t^{j+1,j} \right|^2 \\
&=&\sum_{i=1}^d \left| \sum_{j=0}^{i-1} (Y_t^{j+1,j} - X_t^{j+1,j})  \right|^2 \\
&\leq&\sum_{i=1}^d i \sum_{j=0}^{i-1} |Y_t^{j+1,j} - X_t^{j+1,j}|^2. 
\eeao 

Relation (\ref{NCS-eq:SDsplit_part}) implies
 
$$
 \| Y_t - X_t\|_2^2  \leq \sum_{i=1}^d i\sum_{j=0}^{i-1} (|\wt{Y}_t^{j+1,j}|+ |\wt{X}_t^{j+1,j}|)|\wt{Y}_{t}^{j+1,j} -\wt{X}_{t}^{j+1,j}|.  
$$

Therefore,
\beao
&&\bfE \sup_{0\leq t\leq T} \| Y_t - X_t\|_2^2  \leq \sum_{i=1}^d i \sum_{j=0}^{i-1} \bfE \sup_{0\leq t\leq T} \left((|\wt{Y}_t^{j+1,j}|+ |\wt{X}_t^{j+1,j}|)|\wt{Y}_{t}^{j+1,j} -\wt{X}_{t}^{j+1,j}|\right)\\
&\leq& \sum_{i=1}^d i \sum_{j=0}^{i-1} \sqrt{2\left(\bfE \sup_{0\leq t\leq T}|\wt{Y}_t^{j+1,i}|^2+ \bfE \sup_{0\leq t\leq T}|\wt{X}_t^{i+1,i}|^2\right)}\sqrt{\bfE \sup_{0\leq t\leq T}|\wt{Y}_{t}^{j+1,j} -\wt{X}_{t}^{j+1,j}|^2}\\
&\leq& 4\sqrt{A}\sum_{i=1}^d iC\D^{1/2}\leq C \D^{1/2},
\eeao
where we applied Theorem \ref{NCS-theorem:convergence} and $A=\bfE \sup_{0\leq t\leq T}|\wt{Y}_t^{j+1,i}|^2 \vee \bfE \sup_{0\leq t\leq T}|\wt{X}_t^{i+1,i}|^2$ which is finite by Lemma \ref{NCS-lem:SDuniformMomentBound}.


\subsection*{Final comments and future directions}

Here, we investigated the numerical approximation of systems of SDEs which possess the non-colliding property. The proposed numerical scheme preserves this property. There are however more complicated systems to study;  the interest in this paper was the constant diffusion case. A natural generalization is the following system of SDEs

$$
X_t^{(i)} = X_0^{(i)} + \int_0^t  \left(\sum_{i\neq j}\frac{\gamma_{i,j}}{X_s^{(i)}-X_s^{(j)}} + b^{i}(X_s^{(i)})\right)ds + \sum_{j=1}^d\int_0^t \sigma_{i,j}(X_s^{(i)})dW_s^{(j)}, \quad i=1,\ldots,d,
$$
which arises in various applications such as those in \cite[Section 6]{graczyk_malecki:2014}, \cite{cepa_lepingle_2001}.
We conjecture that one may use the semi discrete method appropriately in a different way. One other possible solution is to use the Lamperti-type transformation as in \cite{neuenkirch_szpruch:2014} to remove the nonlinearity from the diffusion to the drift part of the SDE and then follow the same recipe as the one presented here.


\bibliographystyle{unsrt}\baselineskip12pt 
\bibliography{Non_Colliding_Particle_Systems}

\end{document}